\documentclass[AMA,STIX1COL]{WileyNJD-v2}

\usepackage{amsmath}
\usepackage{amssymb}

\newcommand*{\citen}[1]{%
  \begingroup
    \romannumeral-`\x % remove space at the beginning of \setcitestyle
    \setcitestyle{numbers}%
    [\cite{#1}]%
  \endgroup   
}

    \newcommand{\bfx}{\boldsymbol{x}}

\articletype{RESEARCH ARTICLE}%
\received{}
\revised{}
\accepted{}
\begin{document}

\title{Doubly Degenerate Diffuse Interface Models of Anisotropic Surface Diffusion}

\author[1,2]{Marco Salvalaglio*}
\author[1]{Maximilian Selch}
\author[1,2]{Axel Voigt}
\author[3]{Steven M. Wise}

\authormark{M. SALVALAGLIO \textsc{et al}}

\address[1]{\orgdiv{Institute of Scientific Computing, Department of Mathematics}, \orgname{TU Dresden}, \orgaddress{\state{01062 Dresden}, \country{Germany}}
}
\address[2]{\orgdiv{Dresden Center for Computational Materials Science}, \orgname{TU Dresden}, \orgaddress{\state{01062 Dresden}, \country{Germany}}
}
\address[3]{\orgdiv{Department of Mathematics}, \orgname{The University of Tennessee}, \orgaddress{\state{Knoxville, TN 37996}, \country{USA}}
} 

\corres{*M. Salvalaglio, TU Dresden, 01062 Dresden, Germany. \email{marco.salvalaglio@tu-dresden.de}}

%\maketitle

    \abstract[Summary]{
We extend the doubly degenerate Cahn-Hilliard (DDCH) models for isotropic surface diffusion, which yield
more accurate approximations than classical degenerate Cahn-Hilliard (DCH) models, to the anisotropic case. We consider both weak and strong anisotropies and demonstrate the capabilities of the approach for these cases numerically. The proposed model provides a variational and energy dissipative approach for anisotropic surface diffusion, enabling large scale simulations with material-specific parameters.}


\maketitle

    \section{Introduction}
    \label{sec:introduction}
    
    %Reference to add in future versions
    
    %  Barrett, John W.; Garcke, Harald; Nürnberg, Robert On the stable discretization of strongly anisotropic phase field models with applications to crystal growth. ZAMM Z. Angew. Math. Mech. 93 (2013), no. 10-11, 719–732.

    %Barrett, John W.; Garcke, Harald; Nürnberg, Robert Stable phase field approximations of anisotropic solidification. IMA J. Numer. Anal. 34 (2014), no. 4, 1289–1327.
    
    Surface diffusion is an important transport mechanism in materials science, e.g. in processes such as solid-state dewetting of semiconductors or coarsening of bulk nanoporous metals. In the isotropic setting, surface diffusion is modeled by a geometric evolution equation which relates the normal velocity of a hypersurface in Euclidean space to the surface Laplacian of the mean curvature \cite{Mullins1959}, e.g., $v = \Delta_\Sigma H$. In more realistic anisotropic settings the mean curvature is replaced by a weighted mean curvature, defined as the surface divergence of the Cahn-Hoffmann vector\cite{Cahn1974}, e.g., $H_\gamma = \nabla_\Sigma \cdot D\gamma$, with $\gamma = \gamma(\hat{\mathbf{n}})$ an anisotropic surface energy and $\hat{\mathbf{n}}$ the outward-pointing surface normal. Even though some direct numerical approaches for these equations, in the isotropic and anisotropic setting exist -- see, e.g., \citen{Baenschetal_JCP_2005,Hausseretal_IFB_2005,HausserJSC2007,BarrettNUMMAT2008,Baoetal_JCP_2017,Barrettetal_JCP_2019} -- for most applications in materials science diffuse interface approximations are preferred -- see, e.g., \citen{Cheng2020,Wise2007,Wise2005,Ratz2006,Yeon2006,Torabi2009,Li2009,Banas2009,Jiang2012,Salvalaglio2015a,Bergamaschini2016,Schiedung2017}. These diffuse interface approaches capture the motion of the interface implicitly as the evolution of an iso-surface of a phase-field function. Typically, they are fourth-order nonlinear diffusion equations of Cahn-Hilliard type, whose solutions formally converge to those of their sharp interface counterpart, as the interface thickness tends     to zero\cite{Ratz2006,cahn1996,Gugenberger2008,Dziwniketal_Non_2017}. They require a degenerate mobility function and are termed degenerate Cahn-Hilliard (DCH) equations. In the model proposed in~\citen{Ratz2006}, an additional degeneracy is introduced, following similar ideas as used for the thin film limit in classical phase field models for solidification\cite{Karma1998}. We call such models doubly degenerate Cahn-Hilliard (DDCH) models. This second degeneracy does not alter the asymptotic limit\cite{Ratz2006}, but actually leads to more accurate surface diffusion approximations. See, e.g., the discussion in \citen{Backofen2019}. In fact several simulations for realistic applications in materials science consider this model - see, e.g., \citen{Albani2016,Salvalaglio2015,Salvalaglio2017a,Salvalaglio2017b,Naffouti2017,Geslin2019,Albani2018,Salvalaglio2018,Albani2019,BOL19,BeckAndrewsetal_arXiv_2020,Bergamaschini2020,Salvalaglioetal_PRL_2020} - and claim that their large scale simulations would not be feasible without the additional degeneracy. The drawback of this DDCH model, which is termed \textit{RRV model} in several publications\cite{Gugenberger2008,Geslin2019,BeckAndrewsetal_arXiv_2020}, is that it is non-variational. That is, there is no known free energy that is dissipated along solution trajectories. This makes it harder to prove properties of solutions and derive certain numerical stabilities. In \citen{SalvalaglioDCH} this problem was solved by introducing a new variational DDCH model for surface diffusion, which can be connected with the non-variational DDCH model of \citen{Ratz2006}. This was done for the isotropic case and will here be generalized to the anisotropic case.
    
    The paper is organized as follows: In Sect.~\ref{sec:models} we first review the different isotropic DDCH models, their connection, and then extend the idea to weak as well as strong anisotropies. The latter case requires an additional curvature regularisation, leading to higher-order equations; In Sect.~\ref{sec:numerics} we illustrate the numerical approach which uses a semi-implicit time-integration scheme and adaptive Finite Elements for discretization in space; two- and three-dimensional numerical results are shown in Sect.~\ref{sec:results}, including comparisons between different models and with sharp-interface solutions. Different choices for surface-energy densities are reported, also matching real material properties, and applications are illustrated to show the wide applicability. In Sect.~\ref{sec:conclusions} we draw our conclusions.

    \section{Models} 
    \label{sec:models}
    
    \subsection{Isotropic, variational DDCH model}
    \label{sec:weaklyDDCH}

Suppose that $\Omega\subset\mathbb{R}^d$, $d = 2,3$, is a bounded, open set. Let us consider the free energy 
	\begin{equation}
F[u] = \int_\Omega g_0(u) \left(\frac{1}{\varepsilon}f(u) +\frac{\varepsilon}{2}|\nabla u|^2 \right)d{\bf x} ,
    \label{eqn-reg-energy-g0}
    \end{equation}
with $f= \frac{\omega}{4}u^2(1-u)^2$, $\omega = 72$ a quartic, symmetric double well potential and $g_0$ a singular function of the form
    \[
g_0(u) = \frac{1}{\xi |u|^p |1-u|^p}, \quad p \ge 0, \quad \xi >0.
    \]
$g_0$ can be regularised so that it is defined, continuous, and differentiable for all $u$: 
	\[
g_\alpha(u) = \frac{1}{\sqrt{\xi^2{(u^2(1-u)^2)}^p + \alpha^2\varepsilon^2}}, \quad p\ge 0, \quad \xi >0, \quad \alpha \ge 0.
    \]
The dynamics is described as the $H^{-1}$ gradient flow by 
	\begin{align}
\partial_t u =& \frac{1}{\varepsilon}\nabla\cdot \left(M_0 (u)\nabla w \right),
    \label{eqn1}
    \\
w =& g_0^\prime(u) \left( \frac{\varepsilon}{2} |\nabla u|^2 + \frac{1}{\varepsilon} f(u) \right) + g_0(u) \frac{1}{\varepsilon} f^\prime(u) - \varepsilon \nabla \cdot \left(g_0(u) \nabla u \right),
    \label{eqn2}
 	\end{align}
where $w = \delta_u F$ is the chemical potential with $\delta_u F$ the variational derivative of the free energy $F$ with respect to $u$. $M_0$ is the degenerate mobility function. Its regularised form is
	\begin{equation}
	\label{mob}
M_\alpha(u) = \mu u^2(1-u)^2 + \alpha\varepsilon, \quad \mu = 36, \quad \alpha \ge 0,
    \end{equation}
and $M_0$ is obtained setting $\alpha = 0$. 

In \citen{SalvalaglioDCH} it is shown that Eqs.~\eqref{eqn1} and \eqref{eqn2} formally converge to motion by surface diffusion if $\epsilon \to 0$, $0 \leq p < 2$ and $\xi = 6 \frac{\Gamma(2-p))^2}{\Gamma(4 - 2p)}$ with $\Gamma$ the usual (Bernoulli) Gamma-function. The proposed system, Eqs.~\eqref{eqn1} and \eqref{eqn2}, is a free energy dissipative dynamical system. 

With the asymptotic approximation $\frac{1}{\varepsilon} f(u) \approx \frac{\varepsilon}{2} |\nabla u |^2$ -- which holds when the interface has a hyperbolic tangent profile -- we obtain $w \approx g_0(u) \frac{1}{\varepsilon} f^\prime(u) - \varepsilon g_0(u) \nabla^2 u$, which simplifies Eq.~\eqref{eqn2} and leads to 
	\begin{align}
\partial_t u =& \frac{1}{\varepsilon}\nabla\cdot \left(M_0 (u)\nabla w \right),
    \label{eqn1nv}
    \\
\xi |u|^p |1-u|^p w =& \frac{1}{\varepsilon} f^\prime(u) - \varepsilon \nabla^2 u.
    \label{eqn2nv}
 	\end{align}
This model is strikingly similar to the isotropic version of the model in \citen{Ratz2006}, but with one important caveat. The model in \citen{Ratz2006} corresponds to the choice $p=2$, which is not defined for Eqs.~\eqref{eqn1} -- \eqref{eqn2}\cite{SalvalaglioDCH}.  Interestingly, Eqs.~\eqref{eqn1nv} -- \eqref{eqn2nv} can be defined in a reasonable way for any $p\in [0,\infty)$. The choice of $\xi$ will not be the same as for Eqs.~\eqref{eqn1} -- \eqref{eqn2}, in general. However, for $p=1$ they are identical. See \citen{SalvalaglioDCH} for details.

Numerical solutions indicate a higher accuracy for the variational DDCH model Eqs.~\eqref{eqn1} and \eqref{eqn2} and the non-variational DDCH model Eqs.~\eqref{eqn1nv} and \eqref{eqn2nv} if compared with the classical DCH model ($p = 0$ and $\xi = 1$). Another advantage of the DDCH models is that deviations of $u$ from the pure phase values, 0 and 1, are smaller in a point-wise sense, when compared with the solutions of the classical DCH model. This property can further be elaborated to guarantee, in the non-regularised case $\alpha=0$, that $0 \leq u \leq 1$ for DDCH approximations. So-called positivity preserving schemes of the variational DDCH model Eqs.~\eqref{eqn1} and \eqref{eqn2}, are possible and will be discussed in later publications.

    \subsection{Anisotropic, variational DDCH model}
    \label{sec:anisoDDCH}
    
Instead of a constant surface-energy density $\gamma = 1$, most materials have anisotropic surface-energy density $\gamma = \gamma(\hat{\mathbf{n}})$. In the phase-field context, $\gamma$ is extended by defining  $\hat{\mathbf{n}}^0=-\frac{\nabla u}{|\nabla u|}$.  In other words, $\gamma$ is defined everywhere that $\nabla u\ne {\bf 0}$. We consider the free energy given by 
	\begin{equation}
F[u] = \int_\Omega \gamma(\hat{\mathbf{n}}^0) g_0(u) \left(\frac{1}{\varepsilon}f(u) +\frac{\varepsilon}{2}|\nabla u|^2 \right)d{\bf x}.
    \label{eqn-non-reg-energy}
    \end{equation}
In contrast to classical approach of Kobayashi~\citen{Kobayashi1993}, in which the anisotropic surface-energy density only affects gradient term of the free energy -- see, e.g., \citen{Cheng2020,Wise2007,Wise2005,BarrettNUMMAT2008,Barrett2013,Barrett2014} for stable numerical realisations and even numerical convergence theory~\citen{Cheng2020}-- we here follow the approach of \citen{Torabi2009}, which ensures that the interface thickness is independent of orientation. As in \citen{SalvalaglioDCH}, we include the energy restriction function $g_0$, which is singular at the pure phase values $u = 0,1$. The asymptotic analysis of \citen{SalvalaglioDCH} can be repeated, under the assumption that the diffuse interface has the usual hyperbolic tangent profile, to show that 
$F[u] \approx 
%\int_\Omega \gamma(\hat{\mathbf{n}}_{\Sigma})\delta_\Sigma({\bf x}) \, d{\bf x} =
\int_\Sigma \gamma(\hat{\mathbf{n}}(s)) \, ds,
$
where 
%$\delta_\Sigma$ is the co-dimension one surface delta function associated to 
$\Sigma = \left\{{\bfx}\in\Omega \ \middle| \ u({\bf x}) = 0.5 \right\}$ and $\hat{\bf n}(s)$ is normal for $\Sigma$.

    \subsubsection{Model regularisation}
    
We note that \eqref{eqn-non-reg-energy} is not defined for arbitrary smooth functions $u$. This can be an enormous problem for practical computations. To fix this, we regularise the energy as follows. Set ${\bf p} = \nabla u$. Define the regularised, extended normal
    \[
\hat{\bf n}^\alpha := \frac{-{\bf p}}{\sqrt{|{\bf p}|^2 + \alpha^2\varepsilon^2}},
    \]
and observe that this vector is now (i) only approximately of unit length and (ii) only approximately normal to the level sets of $u$. The regularised free energy is defined as
	\begin{equation}
F_\alpha[u] = \int_\Omega \gamma(\hat{\mathbf{n}}^\alpha) g_\alpha(u) \left(\frac{1}{\varepsilon}f(u) +\frac{\varepsilon}{2}|\nabla u|^2 \right)d{\bf x},
    \label{eqn-reg-energy}
    \end{equation}
which is now well-defined for arbitrary smooth phase-field functions. We note that it is possible, and perhaps even advantageous in some settings, to use two separate  regularisation parameters for $\hat{\bf n}^0$ and $g_0$, but we will avoid this technical discussion here and instead focus on formulation and numerics.

In the computation of the variational derivative of $F_\alpha$, denoted $\delta_u F_\alpha$, we need the following calculation:
    \[
\left[\nabla_{\bf p}\gamma(\hat{\bf n}^\alpha) \right]_i = \frac{\partial\gamma(\hat{\bf n}^\alpha)}{\partial p_i} = \sum_{j = 1}^d P^\alpha_{i,j} \frac{\partial \gamma (\hat{\mathbf{n}}^\alpha)}{\partial n^\alpha_j} = \left[ {\bf P}^\alpha \nabla_{\hat{\bf n}^\alpha} \gamma(\hat{\bf n}^\alpha) \right]_i ,
    \]
where 
    \[
\left[{\bf P}^\alpha \right]_{i,j} = P^\alpha_{i,j} := \frac{\partial n_j^\alpha}{\partial p_i}  = -\frac{ \delta_{i,j} - n^\alpha_i n^\alpha_j }{\sqrt{|{\bf p}|^2+\alpha^2\varepsilon^2}} = -\frac{\left[{\bf I} - {\bf n}^\alpha \otimes {\bf n}^\alpha \right]_{i,j}}{\sqrt{|{\bf p}|^2+\alpha^2\varepsilon^2}} ,
    \]
a regularised projection-like matrix. The regularised doubly degenerate energy dissipative flow now reads
	\begin{align}
\partial_t u & =  \frac{1}{\varepsilon}\nabla\cdot \left(M_\alpha (u)\nabla w \right),
    \label{eqn3}
    \\
w & = \gamma(\hat{\mathbf{n}}^\alpha) \left( g_\alpha^\prime(u) \left( \frac{\varepsilon}{2} |\nabla u|^2 + \frac{1}{\varepsilon} f(u) \right) + g_\alpha(u) \frac{1}{\varepsilon} f^\prime(u) \right) - \nabla \cdot \left( \nabla_{\bf p} \gamma(\hat{\mathbf{n}}^\alpha) g_\alpha(u) \left( \frac{\varepsilon}{2} |\nabla u|^2 + \frac{1}{\varepsilon} f(u) \right) + \gamma(\hat{\mathbf{n}}^\alpha) g_\alpha(u) \varepsilon \nabla u \right),
    \label{eqn4}
 	\end{align}
 where $w = \delta_u F_\alpha$ is the chemical potential. We point out that the chemical potential becomes singular and ill-defined when the regularisation is formally turned off ($\alpha = 0$), in particular, when $\nabla u = {\bf p} = {\bf 0}$,  $u = 0$, or $u = 1$.
 
 	\subsubsection{Strong anisotropy}
 	
Equations~\eqref{eqn3} and \eqref{eqn4} can only be well-posed as long as the graph of $\gamma(\hat{\mathbf{n}})$ is convex, corresponding to so-called weak anisotropies. Strong anisotropies, which also allow for facets, require a corner regularisation. We follow \citen{Torabi2009} and add a diffuse Willmore-type regularisation. The free energy reads 
    \begin{equation}
F_\alpha[u] = \int_\Omega \left\{ \gamma(\hat{\mathbf{n}}^\alpha) g_\alpha(u) \left(\frac{1}{\varepsilon}f(u) +\frac{\varepsilon}{2}|\nabla u|^2 \right)+\frac{\beta}{2\varepsilon } \left(-\epsilon \nabla^2 u +\frac{1}{\varepsilon}f'(u) \right)^2 \right\}  d{\bf x} ,
    \label{eqn-reg-energy-strong}
\end{equation}
with $\beta > 0$. Physically $\beta^{1/2}$ defines a length scale over which corners and edges are smeared out -- see \citen{Torabi2009,Gurtinetal_ARMA_2002}. The considered diffuse interface approximation of the Willmore energy follows from the proposed form in \citen{DeGiorgi1991} by using the asymptotic approximation $\frac{1}{\varepsilon} f(u) \approx \frac{\varepsilon}{2} |\nabla u |^2$. See \citen{Loreti2000,Du_2005,Roeger2006}. The
energy dissipative flow reads
	\begin{align}
\partial_t u=& \frac{1}{\varepsilon}\nabla\cdot \left(M_\alpha (u)\nabla w \right),
    \label{sagf_eqn1}
    \\
w =& \gamma(\hat{\mathbf{n}}^\alpha) \left( g_\alpha^\prime(u) \left( \frac{\varepsilon}{2} |\nabla u|^2 + \frac{1}{\varepsilon} f(u) \right) + g_\alpha(u) \frac{1}{\varepsilon} f^\prime(u) \right)  - \nabla \cdot \left( \nabla_{{\bf p}} \gamma(\hat{\mathbf{n}}^\alpha) g_\alpha(u) \left( \frac{\varepsilon}{2} |\nabla u|^2 + \frac{1}{\varepsilon} f(u) \right) + \gamma(\hat{\mathbf{n}}^\alpha) g_\alpha(u) \varepsilon \nabla u \right) \nonumber \\
&+\beta \left( \frac{1}{\varepsilon^2} f^{\prime\prime}(u)\kappa -\nabla^2 \kappa \right), %+\frac{\beta}{2\varepsilon}g^{\prime}(u) \kappa^2
    \label{sagf_eqn2}
     \\
    \kappa =& -\epsilon \nabla^2 u +\frac{1}{\varepsilon} f'(u).
    \label{sagf_eqn3}
	\end{align}
We remark that using a \emph{low-order} mobility function -- such as $M_0(u) = Cu(1-u)$ which  was employed in \citen{Torabi2009} -- solutions would not actually converge to motion by anisotropic surface diffusion, as $\epsilon \to 0$. The reasons are laid out in  \citen{dai2014,Leeetal_APL_2015,Leeetal_SIAMJAM_2016}. However, the asymptotic analysis with $M_0(u)$ as in Eq.~\eqref{mob} leads to the desired results for both models \eqref{eqn3} - \eqref{eqn4} and \eqref{sagf_eqn2} - \eqref{sagf_eqn3}, provided $g_0 = 1$, as shown in \citen{Ratz2006,Gugenberger2008,Voigt2016}. Additionally, combining these results with the results of \citen{SalvalaglioDCH} will give the same asymptotic limits as $\varepsilon\to 0$ for the doubly degenerate cases with $g_0(u)$. The limiting model reads $v = \Delta_\Sigma (H_\gamma + \beta (\Delta_\Sigma H + \frac{1}{2} H^3 - 2HK))$, with Gaussian curvature $K$ -- see, e.g., \citen{Gurtinetal_ARMA_2002}.

Also for the anisotropic cases the asymptotic approximation $\frac{1}{\varepsilon} f(u) \approx \frac{\varepsilon}{2} |\nabla u |^2$ can be used to simplify the equations and, in turn, their numerical integration. However, a one-to-one correspondence with the non-variational models proposed in \citen{Ratz2006} is not possible. The corresponding version, including the surface energy density as in \citen{Torabi2009}, reads 
    \begin{align}
\partial_t u =& \frac{1}{\varepsilon}\nabla\cdot \left(M_0 (u)\nabla w \right),
    \label{nv-eqn3}
    \\
\xi |u|^p |1-u|^p w =& -\varepsilon \nabla \left( \gamma(\hat{\mathbf{n}}^0) \nabla u \right) + \frac{1}{\varepsilon} \gamma(\hat{\mathbf{n}}^0) f'(u) - \varepsilon \nabla \cdot \left( |\nabla u|^2 \nabla_{\mathbf{p}} \gamma(\hat{\mathbf{n}}^0) \right) + \beta \left( \frac{1}{\varepsilon^2} f^{\prime\prime}(u)\kappa -\nabla^2 \kappa \right),
    \label{nv_eqn4}
    \\
    \kappa =& -\epsilon \Delta u +\frac{1}{\epsilon} f'(u),
    \label{nv_eqn5}
	\end{align}
with $p =1$, $\xi = 6$ or $p = 2$, $\xi = 30$, and $\beta = 0$ for weak and $\beta > 0$ for strong anisotropy. Also these models formally converge asymptotically to the expected anisotropic surface diffusion models as $\varepsilon\to 0$ and have been frequently used in applications - see, e.g., \citen{Salvalaglio2015a,Bergamaschini2016,BOL19}. To be practically useful also Eqs.~\eqref{nv-eqn3} -- \eqref{nv_eqn5} require an appropriate regularisation \cite{Ratz2006}. In the following, we will refer to Eqs.~\eqref{nv-eqn3} -- \eqref{nv_eqn5} as the \textit{RRV model}.

    \section{Numerics}
    \label{sec:numerics}

Our goal is to demonstrate the advantageous solution properties of the DDCH models, as well as to show that the variational DDCH models are suitable for simulating the formation of complex faceted morphologies. We will in the following consider the case $p = 1$ and $\xi = 6$ for the variational model \eqref{sagf_eqn1} -- \eqref{sagf_eqn3}. The classical DCH model can be straightforwardly considered by setting $p=0$ and $\xi=1$. (Thus, we can take $g_\alpha \equiv 1$). For the time-integration of the \textit{RRV model} we refer to \citen{Ratz2006,Torabi2009}. We employ the regularised functions $\hat{\bf n}^\alpha$, $g_\alpha$, and $M_\alpha$ using the value $\alpha=10^{-6}$. 
%For the derivative of $\gamma$, it is convenient to consider the following identity
%    \begin{equation}
%    \partial_{\nabla u} \gamma(\hat{\mathbf{n}}) = \frac{1}{|\nabla u|}\mathbf{P} \nabla_{\hat{\mathbf{n}}} \gamma(\hat{\mathbf{n}}),
%    \label{eq:gradientGamma}
%    \end{equation}
%with $\mathbf{P} = \mathbf{I} - \hat{\mathbf{n}} \otimes \hat{\mathbf{n}}$ and $\nabla_{\hat{\mathbf{n}}} = (\frac{\partial}{\partial \hat{{n}}_x},\frac{\partial}{\partial \hat{{n}}_y},\frac{\partial}{\partial \hat{{n}}_z})$. 

We solve the models using adaptive finite elements, implemented in the software package AMDiS\cite{Vey2007,Witkowski2015}. The time-integration scheme is semi-implicit and follows the scheme employed for the isotropic case\cite{SalvalaglioDCH}, which reads,
 \begin{align}
\frac{u^{n+1}}{\tau_n} - \frac{1}{\varepsilon} \nabla \cdot \left( M_\alpha(u^{n}) \nabla w^{n+1}\right)&= \frac{u^{n}}{\tau_n},
     \label{eq:system_iso1}
   \\
w^{n+1}+\varepsilon\nabla \cdot \left[ g_\alpha(u^{n}) \nabla u^{n+1}\right]-\frac{1}{\varepsilon}\left[r(u^{n})+ s(u^{n})\right] u^{n+1}
+\frac{\varepsilon}{2} g_\alpha'(u^{n})\nabla u^{n} \cdot \nabla u^{n+1}&=
q(u^n),
    \label{eq:system_iso2}
    \end{align}
where
    \begin{align*}
q(u^n)&=\frac{1}{\varepsilon}g_\alpha(u^{n})f'(u^{n})+\frac{1}{\varepsilon}g'_\alpha(u^{n})f(u^{n})-\frac{1}{\varepsilon}\left[r(u^{n})+s(u^{n})\right]u^{n},
    \\
r(u^{n})&=[g'_\alpha(u^{n})f'(u^{n})+g_\alpha(u^{n})f''(u^{n})],
    \\
s(u^{n})&=[g'_\alpha(u^{n})f'(u^{n})+g''_\alpha(u^{n})f(u^{n})], 
    \end{align*}
account for the linearizations of $g_\alpha(u^{n+1})f'(u^{n+1})$ and $g'_\alpha(u^{n+1})f(u^{n+1})$ around $u^{n}$. The integer $n$ is the time step index, $\tau_n$ the time stepsize at step $n$, and $u^0$ the initial condition. 

To employ this scheme for the anisotropic cases, $\gamma(\hat{\mathbf{n}}^\alpha)$ has to be included. The scheme for the variational model \eqref{sagf_eqn1} -- \eqref{sagf_eqn3} reads
    \begin{align}
\frac{u^{n+1}}{\tau_n} - \frac{1}{\varepsilon}\nabla \cdot \left[ M_\alpha(u^{n}) \nabla w^{n+1}\right] &= \frac{u^{n}}{\tau_n},   \label{eq:system_full_strong3a} 
    \\
w^{n+1}+\varepsilon \nabla \cdot \left( \gamma_\alpha^n g_\alpha(u^n) \nabla u^{n+1} \right) - \gamma_\alpha^n \frac{1}{\varepsilon} \left[ r(\varphi^n) + s(\varphi^n) \right] u^{n+1} \nonumber 
    \\ 
+ \frac{\varepsilon}{2} \gamma_\alpha^n g_\alpha'(u^{n})\nabla u^{n} \cdot \nabla u^{n+1}
    + \beta \left( - \frac{1}{\varepsilon^2} 
    %g(u^{n})
    f^{\prime\prime}(u^n) \kappa^{n+1} + 
    %g(u^{n})
    \nabla^2 \kappa^{n+1}
    %-\frac{\beta}{2\varepsilon}g^{\prime}(u^{n}) \kappa^n
    \right)  & =  \gamma_\alpha^n q(u^n) + A(u^n),     \label{eq:system_full_strong3b} 
    \\
\kappa^{n+1} + \varepsilon \nabla^2 u^{n+1} -\frac{1}{\varepsilon}f^{\prime\prime}(u^n)u^{n+1}&= \frac{1}{\varepsilon} f^\prime(u^n)- \frac{1}{\varepsilon}f^{\prime\prime}(u^n)u^n,
    \label{eq:system_full_strong3c} 
    \end{align}
where
    \[
\gamma_\alpha^n := \gamma\left(\hat{\mathbf{n}}^{\alpha,n}\right), \quad \hat{\mathbf{n}}^{\alpha,n} :=- \frac{{\bf p}^n}{\sqrt{|{\bf p}^n|^2+\alpha^2\varepsilon^2}}, \quad {\bf p}^n := \nabla u^n,
    \]
    \begin{equation}
A(u^n) :=- \nabla \cdot \left( g_\alpha(u^n)\mathbf{P}^{\alpha,n} \nabla_{\hat{\mathbf{n}}^\alpha} \gamma(\hat{\mathbf{n}}^{\alpha,n}) \left( \frac{\varepsilon}{2} |\nabla u^n|^2 + \frac{1}{\varepsilon} f(u^n) \right) \right) , \quad \mbox{and} \quad  \left[{\bf P}^{\alpha,n} \right]_{i,j} = P^{\alpha,n}_{i,j} := -\frac{ \delta_{i,j} - n^{\alpha,n}_i n^{\alpha,n}_j }{\sqrt{|{\bf p}^n|^2+\alpha^2\varepsilon^2}} .
   \label{eq:Aun}
    \end{equation}
%The model above has a singular term, owing to the derivative of $\gamma$. See eqn.~\eqref{eq:gradientGamma}. To avoid dividing by zero, we regularise the gradient term $|\nabla u^n| \rightarrow |\nabla u^n| + \alpha \epsilon$. 

For weak anisotropies we set $\beta = 0$ and Eq.~\eqref{eq:system_full_strong3c} is not needed. The system is discretized in space by piecewise linear finite elements. We consider simplicial meshes which are adaptively refined by bisection to ensure a resolution ranging between $5 h \sim \varepsilon$ and $10 h \sim \varepsilon$, with mesh size $h$, within the diffuse interface. The resulting linear system is solved by an iterative solver exploiting the biconjugate gradient stabilized method (BiCGstab(l))\cite{VDV1992}, with a block Jacobi preconditioner (bjacobi) applied to the blocks resulting from an element-wise domain decomposition, with a local sparse direct solver UMFPACK\cite{Davis2004}.

    \begin{remark}
Since we employ the regularized models in our numerical experiments, the positivity of $u$, namely, $0<u<1$, is neither expected to hold, nor required to hold, in the numerical results. However, as the regularization parameter is decreased, the degree of the overshoots -- that is, $u$ going above 1 and/or below 0 -- are reduced. We have confirmed this in numerical tests not reported here. A comparison of overshoots in the various models is graphically illustrated in Figure~\ref{fig:figure2}. The design of a theoretically positivity preserving and energy dissipating numerical scheme for the present anisotropic variational model would be a difficult undertaking, in our experience. However, it is possible to develop such a scheme for the isotropic variational DDCH model, and we plan to give the details of that in a future paper.
    \end{remark}

\section{Numerical Results}
\label{sec:results}

Numerical results are chosen to demonstrate the advantageous solution properties and possibilities of the variational DDCH model for applications in materials science. We compare the variational DDCH model with the \textit{RRV model} and the classical DCH model for various anisotropies in 2D. The chosen examples in 3D show faceting with material specific anisotropies and large scale coarsening simulations in nano-porous materials. 

    \begin{figure}[h]
    \center
\includegraphics*[width = \textwidth ]{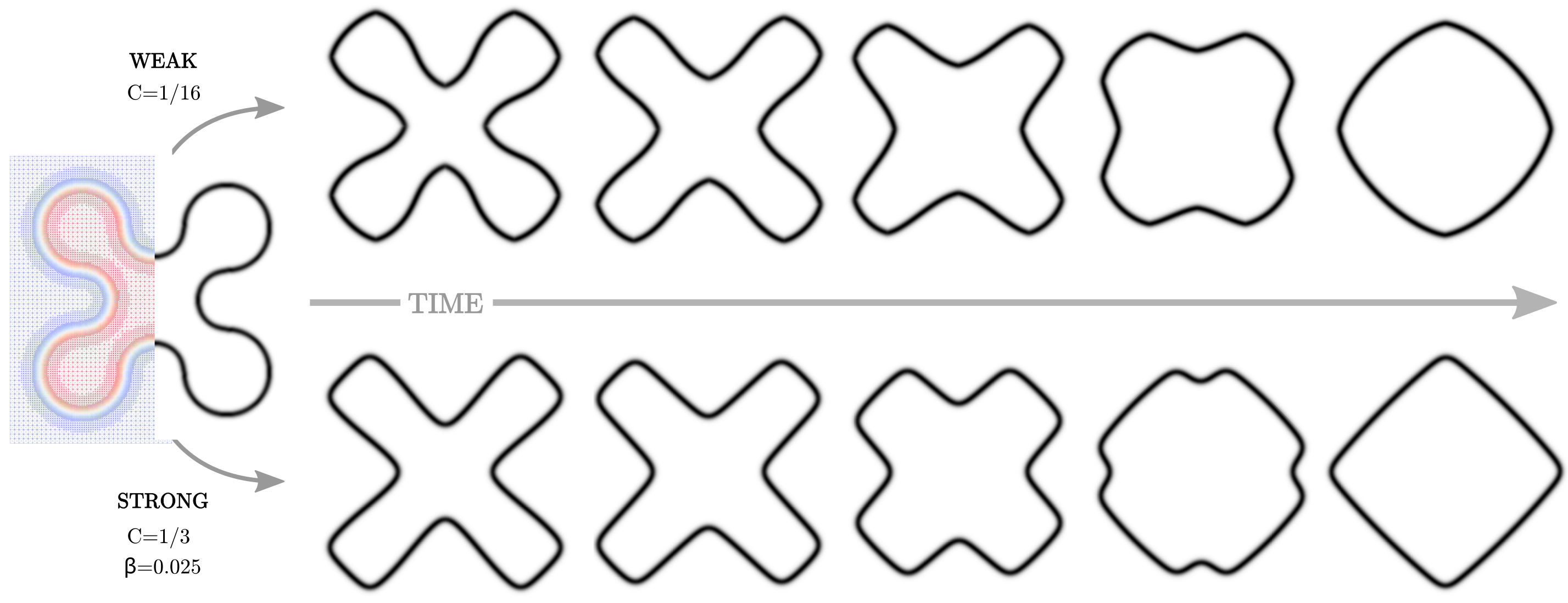} 
\caption{Evolution obtained by the anisotropic variational DDCH model ($\varepsilon=0.2$). The shapes are obtained as white-to-black color maps for $|\nabla u|^2$. Different shapes (left to right) correspond to time $0\tau$, $10^{2}\tau$, $5\cdot 10^{2}\tau$, $15\cdot 10^{2}\tau$, $30\cdot 10^{2}\tau$, $60\cdot 10^{2}\tau$ with $\tau=2\cdot 10^{-5}$ the timestep for weak anisotropy and $\tau=1\cdot 10^{-5}$ the timestep for strong anisotropy. The first shape on the left corresponds to the initial condition; it shows also the adaptive computational mesh and the bounds of the computational domain with size $2 \times 2$. The last shapes on the right are the reached equilibrium shapes. Anisotropic surface energy is set as in Eq.~\eqref{eq:gammaForm} with parameters as specified in the figure. The solutions are independent to further mesh refinement.}
    \label{fig:figure1}
    \end{figure}

\begin{figure}[h]
    \center
\includegraphics*[width = 0.8\textwidth ]{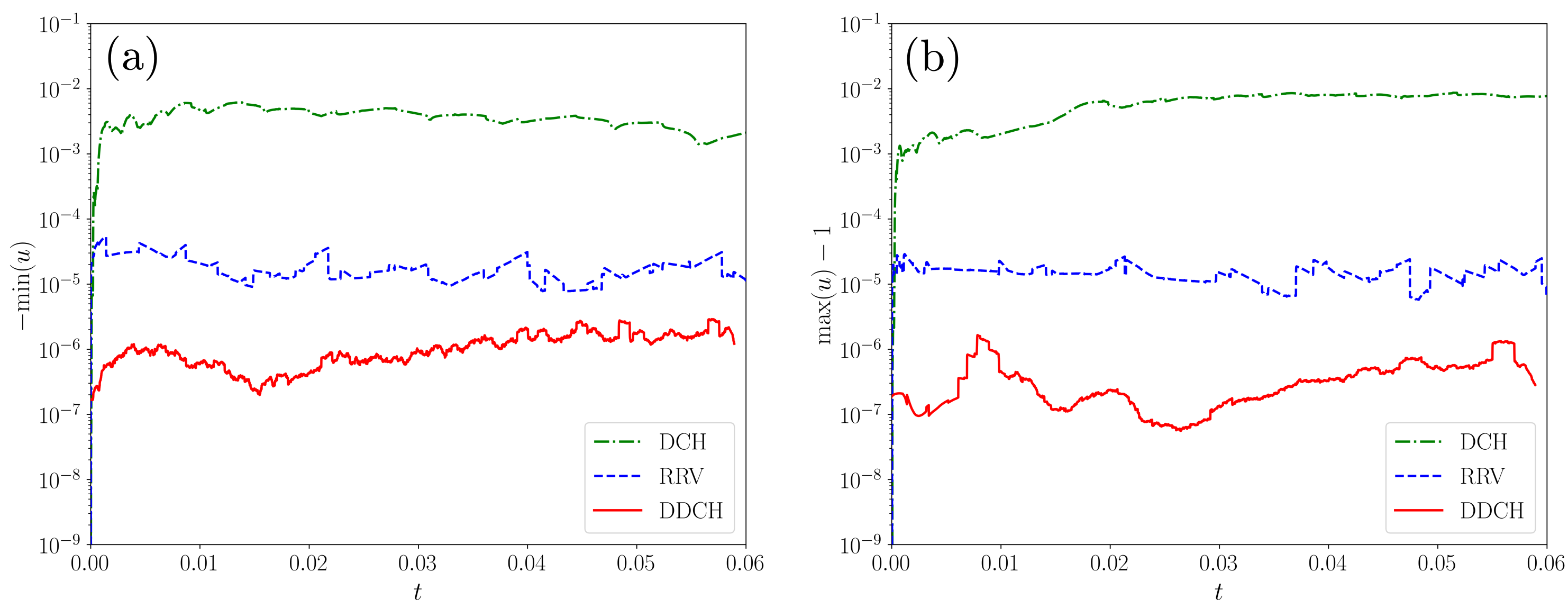} 
\caption{The minimum and maximum of $u$ for the simulation in Fig.~\ref{fig:figure1} (weak anisotropy) is shown by means of $-\text{min(u)}$ (a) and $\text{max(u)}-1$ (b), respectively, for: i) the standard degenerate Cahn-Hillidard model (DCH, green dashed-dotted line), ii) the RRV model (dashed blue line), iii) the variational doubly degenerate Cahn-Hilliard model Eqs.~\eqref{eqn3} and \eqref{eqn4} (DDCH, red solid line).}
\label{fig:figure2}
\end{figure}

In 2D we consider a simple, four-fold anisotropy
    \[
\gamma_4(\hat{\bf n}) =  1+ C\left(4\left(n_x^4 + n_y^4\right) -3 \right),
    \]
which in polar form is
    \begin{equation}
\gamma_4(\hat{\mathbf{n}}) = \tilde\gamma_4(\theta)=1+C\cos(4\theta), \quad \theta=\text{atan}\left(\frac{{n}_y}{{n}_x}\right) .
    \label{eq:gammaForm}
    \end{equation}
With this choice, $\gamma(\hat{\mathbf{n}})$ is convex (weakly anisotropy) for values $C<\frac{1}{15}$ and non-convex (strongly anisotropy) for $C>\frac{1}{15}$. The evolution obtained with Eqs.~\eqref{eqn3} - \eqref{eqn4} or Eqs.~\eqref{sagf_eqn2} - \eqref{sagf_eqn3}, depending on $C$, is illustrated in Fig.~\ref{fig:figure1}. The initial setting is a 2D shape with varying positive and negative curvatures as in \citen{SalvalaglioDCH}. The final shape, the obtained equilibrium configuration, is the approximation of the corresponding Wulff shape, with facets and rounded corners for the non-convex (strongly anisotropic) free energy. For both cases the evolution can also be obtained with the RRV and DCH models. (Recall that the DCH model is obtained from the DDCH model using $g_\alpha \equiv 1$.) Visually there is no difference. However, as already discussed for the isotropic case, the variational DDCH model has better ``positivity'' properties. Note, however, that since we are using a regularised model, and since our numerical method is not designed for positivity, the numerical approximations of the  DDCH model do not preserve the positivity property $0<u<1$ precisely. However the overshoot values are very small for the DDCH and RRV models. The DDCH model performs the best, as the overshoots are closer to the pure phase values $u = 0$ and $u = 1$.

\begin{figure}[h]
    \center
\includegraphics*[width = \textwidth ]{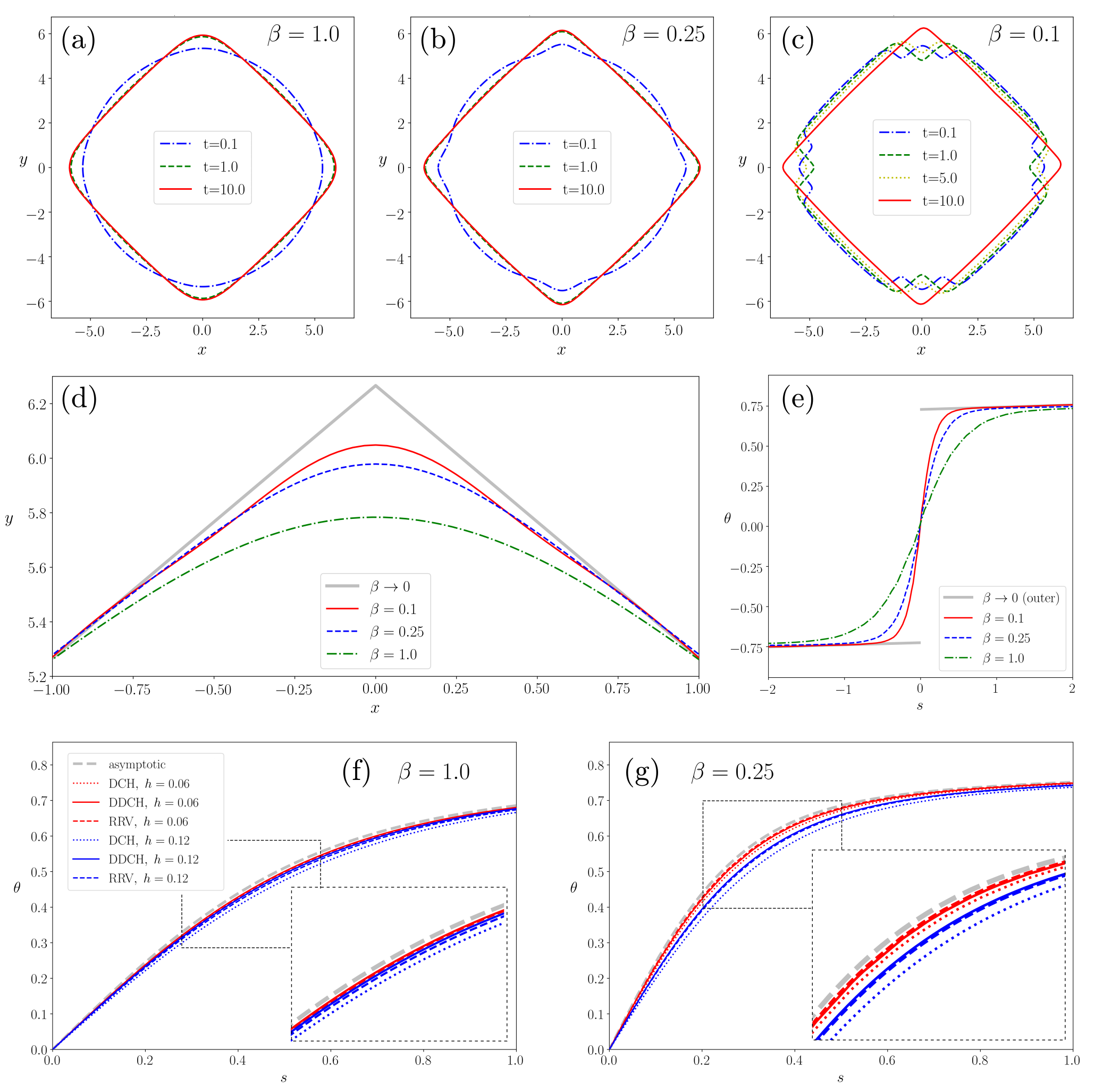} 
\caption{Analysis of faceting and corner rounding. Evolutions from a circle with radius 5 ($\varepsilon=0.6$) to equilibrium with anisotropy as in Eq.~\eqref{eq:gammaForm} with $C = \frac{1}{2}$: (a) $\beta=1.0$, (b) $\beta=0.25$ and (c) $\beta=0.1$. The solutions are independent to further mesh refinement. Comparison of corner rounding in corresponding Wulff shape of limiting sharp interface problem if $\varepsilon \to 0$ for different $\beta$. Shown is the upper corner in (x,y) plot (d), and orientation $\theta$ as function of the arclength $s$ (e). Comparison of $\theta(s)$ close to the corner obtained by the variational DDCH, RRV and DCH models for two spatial discretizations for $\beta=1.0$ (f) and $\beta=0.25$ (g). The legend in (f) also holds for (g).}
    \label{fig:figure3}
\end{figure}
    
Clear evidence of this aspect in the anisotropic case is reported in Fig.~\ref{fig:figure2}. It shows the deviation from the nominal values of $u$, namely 0 and 1, away from the interface through $-\text{min}(u)$ and $\text{max}(u)-1$, respectively. The variational DDCH model shows values, which are between one and two orders of magnitudes closer to the nominal values (0 in both plots) with respect to the \textit{RRV model} and three to four orders of magnitude closer with respect to the DCH model. It is worth mentioning that these values depend on the specific simulation set up. However, we expect the relative difference between the models to hold true generally.

The better approximation properties of the variational DDCH and RRV models, if compared with the classical DCH model in the isotropic case, have already been shown in \citen{SalvalaglioDCH,Backofen2019}, respectively. These properties result in the possibility to use larger $\epsilon$ and thus lower grid resolution, to reach the same accuracy, which enables the large-scale simulations in various applications. In the anisotropic setting the differences between the variational DDCH and RRV models and the classical DCH model become even more severe and might lead to divergence of the DCH model for strong anisotropies and desired error tolerances. To demonstrate this we compare the reached equilibrium shape with the desired Wulff shape of the sharp interface problem. We closely follow convergence studies considered in \citen{Hausseretal_IFB_2005} for a sharp interface model for strongly anisotropic surface diffusion. Fig. \ref{fig:figure3}(a) - (c) show the evolution from a circle with $\text{radius} = 5$ towards the equilibrium shape for different $\beta$. We consider $\varepsilon = 0.6$ and $\gamma(\hat{\mathbf{n}})$ as in Eq.~\eqref{eq:gammaForm} with $C=\frac{1}{2}$. Focusing on the rounded upper corner of the desired Wulff shape allows to construct an asymptotic solution for the limiting sharp interface problem, see \citen{Spencer2004} for details. 
Briefly, the idea is to take the sharp corner equilibrium shape ($\beta = 0$) as the outer solution and to derive an inner solution for the equilibrium shape near the corner as an expansion in $\beta^{1/2}$. Let $s$ be the arclength with $s = 0$ at the corner of the outer solution and $\Theta(S) = \theta(s / \beta^{1/2})$ the rescaled local orientation. Expanding $\Theta(S) = \Theta_0(S) + \beta^{1/2} \Theta_1(S) + \cdots$, one obtains the lowest order term $\Theta_0(S)$ by inverting 
\begin{equation}
    S(\Theta_0)=\int_0^{\Theta_0} \frac{1}{\sqrt{2Q(\Theta')}}d\Theta' \qquad \text{with} \qquad Q(\Theta_0)= \gamma(\Theta_0) + A \cos(\Theta_0)
\end{equation}
with $A=-\gamma(\theta_c)\cos(\theta_c)+\gamma'(\theta_c)\sin(\theta_c)$ and $\theta_c=\pm \frac{\pi}{4}$ the corner orientation of the outer solution. The composite solution in the neighborhood of a corner with inner solution $\theta_{\rm inner}$ is then given by
$\theta(s)=\theta_{\rm inner}(s)+\theta_{\rm outer}(s)-\theta_{\rm match}$, with $\theta_{\rm match} = \pm \frac{\pi}{4}$, with the sign depending on the sign of $s$. Fig. \ref{fig:figure3}(d),(e) show the solution for various $\beta$ in an $(x,y)$- and $(s,\theta)$-plot, respectively. The comparison of our numerical solutions of the variational DDCH, RRV and DCH models for fixed $\varepsilon$ with these asymptotic solutions is shown in Fig.~\ref{fig:figure3}(f),(g), for $\beta=1$ and $\beta=0.25$, respectively. Two spatial discretizations are shown corresponding to 5 ($h=0.12$ and 10 ($h=0.06$) grid points across the diffuse interface, demonstrating the good approximation properties of the double degenerate models already for moderate values of $\varepsilon$.

\begin{figure}[h]
    \center
\includegraphics*[width = 0.9\textwidth ]{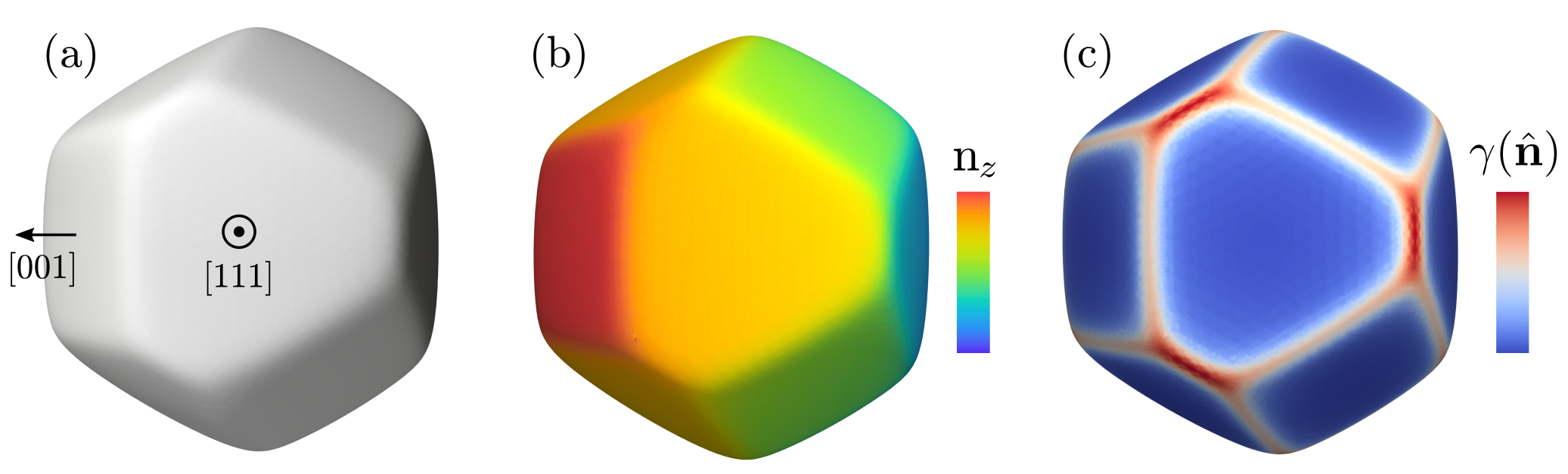} 
\caption{Equilibrium crystal shape obtained with the variational DDCH model. The initial configuration is a sphere with radius 1 ($\varepsilon=0.2$). $\gamma(\hat{\mathbf{n}})$ as in Eq.~\eqref{eq:gammaCGD} with ${\mathbf{m}}_i$ corresponding to $\langle 001 \rangle$ and $\langle 111 \rangle$ directions, $\gamma_0=1$, $\rho_i=0.3$, $\omega_i=10$, and $\beta=0.005$. (a) Isosurface $u = \frac{1}{2}$. (b) z-component of $\hat{\mathbf{n}}$, highlighting the formation of preferential orientations. (c) $\gamma(\hat{\mathbf{n}})$.}
    \label{fig:figure4}
    \end{figure}

For examples in 3D with strong anisotropy we consider $\gamma(\hat{\mathbf{n}})$ as proposed in \citen{Salvalaglio2015a}: 
\begin{equation}
 \gamma(\hat{\mathbf{n}})=\gamma_0 \left[1 - \sum_{i}^N \rho_i \left(\hat{\mathbf{n}} \cdot\mathbf{m}_i\right)^{\omega_i} \Theta \left(\hat{\mathbf{n}}\cdot\mathbf{m}_i \right) \right],
 \label{eq:gammaCGD}
\end{equation}
with $\mathbf{m}_i$ the normals to the surface corresponding to minima of $\gamma(\hat{\mathbf{n}})$, $\rho_i$ the depth of energy minima, $\omega_i$ controlling the extension of energy well around $\mathbf{m}_i$ and $N$ the total number of different minima. This particular form provides the possibility to specify anisotropies for a huge class of materials. Notice that with $\omega=4$, $\rho_i=\bar{\rho}$, $\gamma_0=1$, $\mathbf{m}_i=\pm \mathbf{e}_i$ this form for $\gamma(\hat{\mathbf{n}})$ corresponds to another well-known four-fold surface-energy function,
    \begin{equation}
\gamma_4(\hat{\mathbf{n}})= 1 - \bar{\rho}(n_x^4+n_y^4+n_z^4).
    \label{eq:simpleg4}
    \end{equation} 
An example of faceting of a sphere in a strong anisotropy regime is shown in Fig.~\ref{fig:figure4}. Therein, a surface energy density with minima along $\{111\}$ and $\{100\}$ directions has been considered (see parameters in the caption). The shape corresponding to the isosurface $u=\frac{1}{2}$ is shown (panel (a)), along with the z-component of $\hat{\mathbf{n}}$ and $\gamma(\hat{\mathbf{n}})$ at the same isosurface (panel (b) and (c), respectively), highlighting the facets along with the rounded edges and corners. 

\begin{figure}[t]
    \center
\includegraphics*[width = \textwidth ]{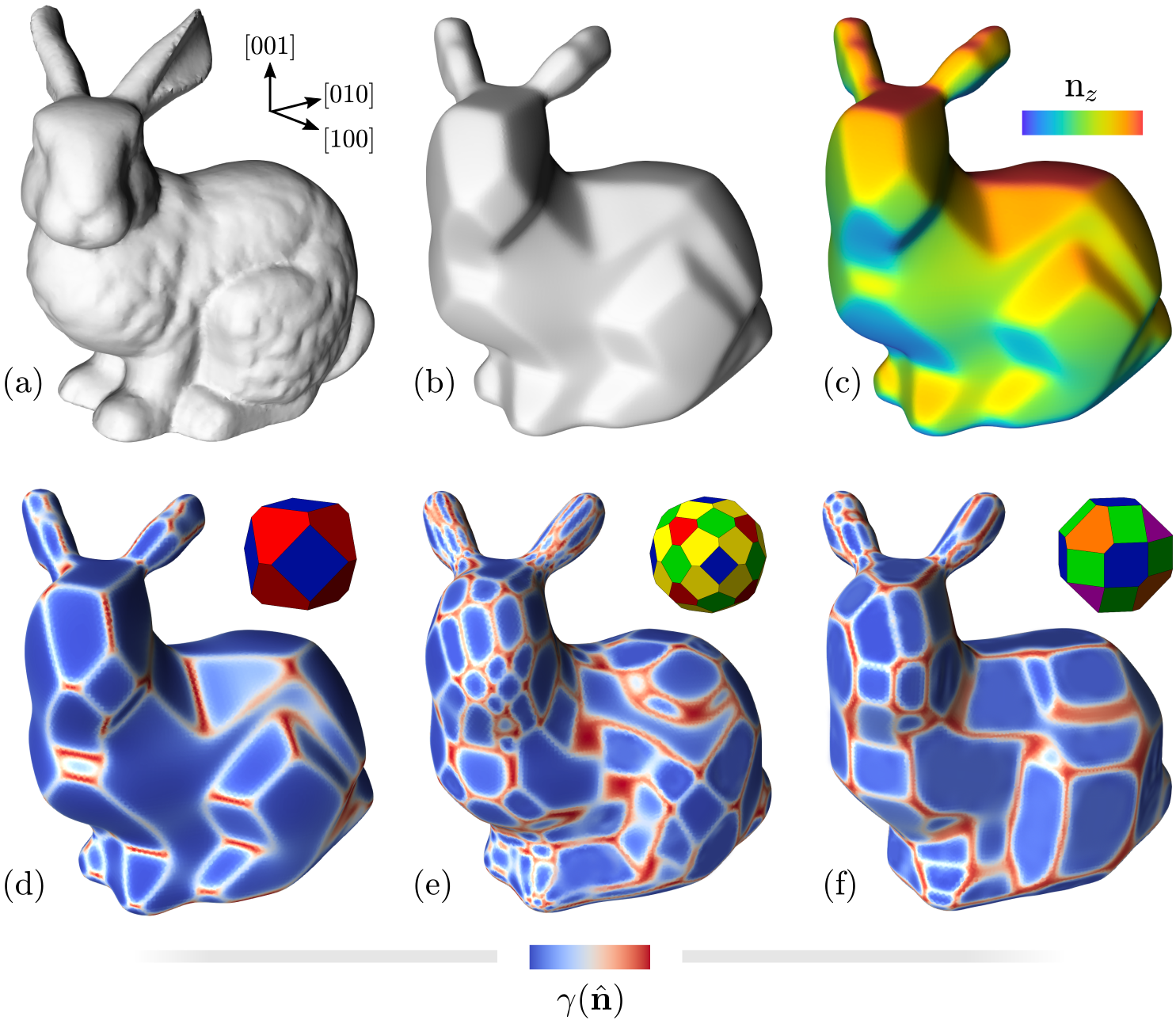} 
\caption{Faceting of the Stanford Bunny\cite{BUNNY} by the DDCH model. The computational domain with size $1 \times 1 \times 1$ ($\varepsilon=0.1$). (a) Initial condition (isosurface $u= \frac{1}{2}$); (b) Faceted geometry after several time steps with anisotropy as in Fig.~\ref{fig:figure4} with $\beta=2\cdot 10^{-4}$; (c) z-component of $\hat{\mathbf{n}}$; (d) $\gamma(\hat{\mathbf{n}})$; (e) Faceted geometry, illustrated by $\gamma(\hat{\mathbf{n}})$, obtained with a surface-energy anisotropy parametrization as for Si crystals with $\beta=3\cdot 10^{-4}$; (f) as for GaAs crystals with $\beta=3\cdot 10^{-4}$. The Wulff shape (without edge and corner rounding) corresponding to anisotropies of panel (d)--(f) is also reported therein for comparison.}
    \label{fig:figure5}
    \end{figure}

The faceting of a more complex shape is shown in Fig.~\ref{fig:figure5}. We considered the \textit{Stanford bunny}\cite{BUNNY} (Fig.~\ref{fig:figure5}(a)). The same surface energy anisotropy as in Fig.~\ref{fig:figure4} with $\beta=2\cdot 10^{-4}$ is considered in Fig.~\ref{fig:figure5}(b)-(d), illustrating the $u=\frac{1}{2}$ isosurface, the z-component of $\hat{\mathbf{n}}$ and $\gamma(\hat{\mathbf{n}})$ after several time steps, respectively. Fig.~\ref{fig:figure5}(e)-(f) show the corresponding solution for different anisotropies. The parameters are set to reproduce the main facets of Si and GaAs crystals, respectively. The parameters considered for Si crystals: $\mathbf{m}_i$=$\{\langle001\rangle,\langle 111\rangle,\langle110\rangle,\langle113\rangle \}$, $2\omega_{\langle001\rangle}$=$2\omega_{\langle111\rangle}$=$2\omega_{\langle110\rangle}$=$\omega_{\langle113\rangle}$=$100$; $\rho_{\langle001\rangle}=0.08$, $\rho_{\langle111\rangle}=0.07$, $\rho_{\langle110\rangle}=0.07$ $\rho_{\langle113\rangle}=0.075$ (qualitatively reproducing energetics as in \citen{Gai1999}), and for GaAs crystals: $\mathbf{m}_i$=$\{\langle001\rangle,\langle 111\rangle,\langle110\rangle\}$, $\omega_{\langle001\rangle}$=$\omega_{\langle111\rangle}$=$\omega_{\langle110\rangle}$=$50$, $\rho_{\langle001\rangle}=0.1$, $\rho_{\langle111{\rm A}\rangle}=0.06$, $\rho_{\langle111{\rm B}\rangle}=0.12$ $\rho_{\langle110\rangle}=0.12$ (qualitatively reproducing energetics as in \citen{Moll1996}, neglecting $\{113\}$ facets). See also \citen{Salvalaglio2017a,Salvalaglio2018}. 

The nominal equilibrium crystal shapes, without corner rounding induced by the Willmore regularisation, are reported in the upper right corners of Fig.~\ref{fig:figure5}(d)-(f) for comparison. Notice that the shapes obtained in Fig.~\ref{fig:figure5} are not stationary. They show the solution at one time instance of the evolution, depicted to visualize the faceting of the initial condition. To comment on the positivity preserving property of the variational DDCH model, we have computed the average values for $-\text{min}(u)$ and $\text{max}(u)-1$ over all simulations in Fig.~\ref{fig:figure5}. We obtain the respective values $(1.1 \pm 0.4) \cdot 10^{-5}$ and $(0.8 \pm 0.3) \cdot 10^{-5}$, which are a bit larger than those exhibited in Fig~\ref{fig:figure2} but still within the range of the \textit{RRV model} in Fig~\ref{fig:figure2}. This deviation is expected to be larger in the 3D setting. 

\begin{figure}[h]
    \center
\includegraphics*[width = 0.32\textwidth ]{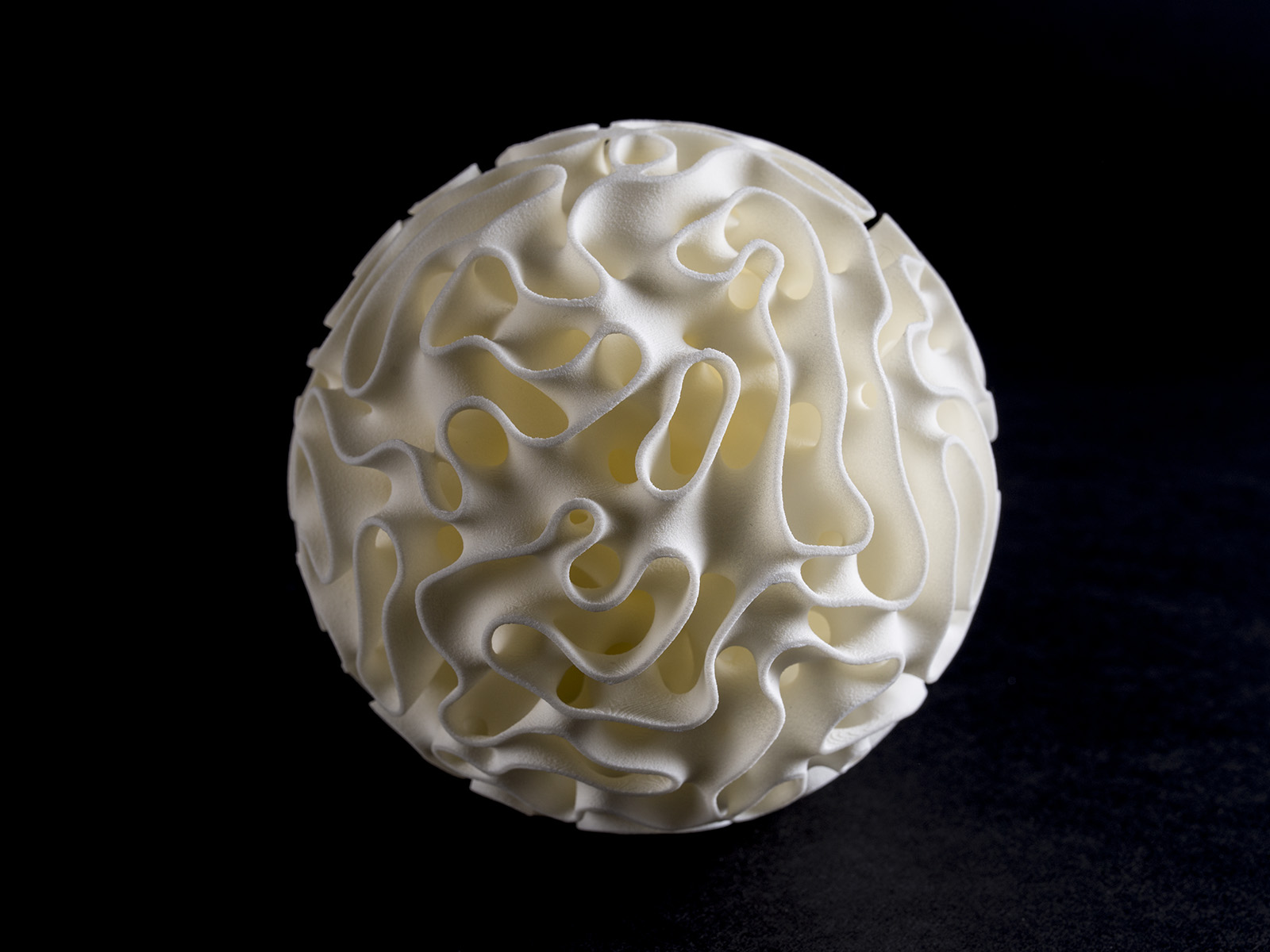} 
\includegraphics*[width = 0.32\textwidth ]{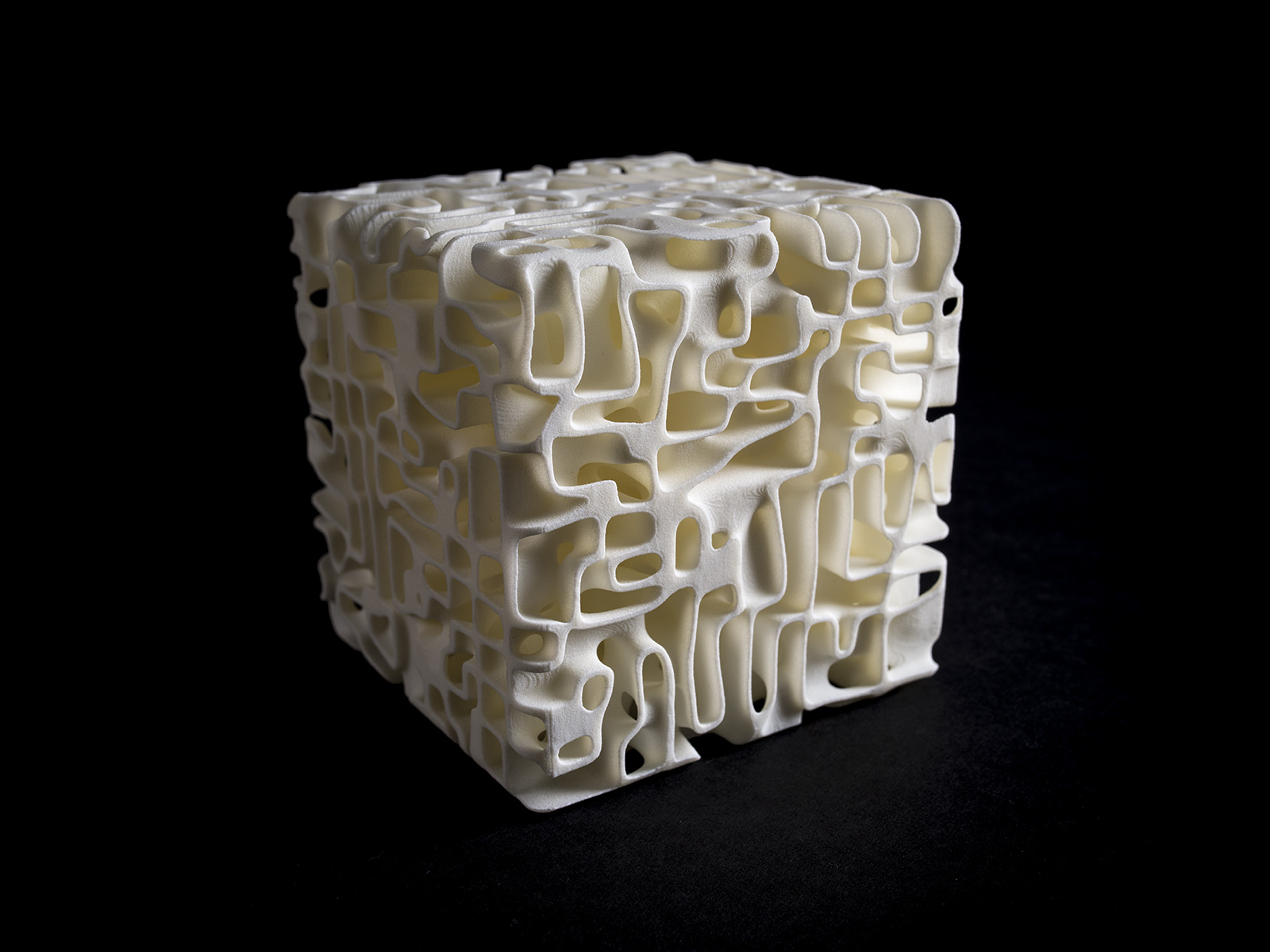} 
\includegraphics*[width = 0.32\textwidth ]{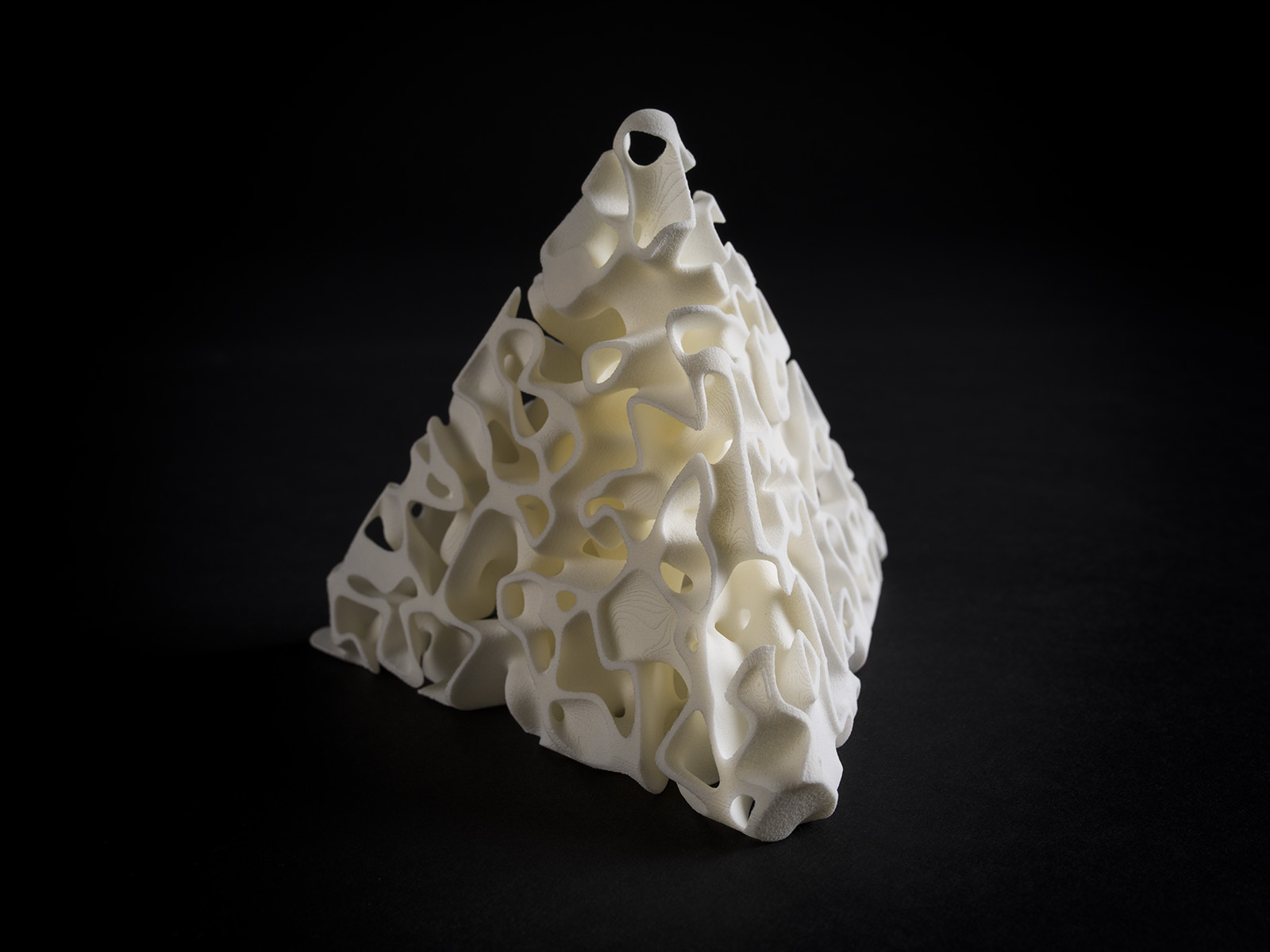} 
\caption{
Objects are part of the Digital Archive of Mathematical Models (DAMM)~\cite{damm}, computed using AMDiS\cite{Vey2007,Witkowski2015} by R. Backofen and F. Stenger, 3D print by materialise.com, photo by D. Lordick. Surface-energy anisotropies are set by Eq.~\eqref{eq:gammaCGD}. Left panel: isotropic setting. Center and Right panels: cubic faceting with $\mathbf{m}_i=\pm \mathbf{e}_i$ and a tetrahedral faceting with 
    $\mathbf{m}_1=[\bar{1}\bar{1}\bar{1}]$,
$\mathbf{m}_2=[11\bar{1}]$,
$\mathbf{m}_3=[1\bar{1}1]$,
$\mathbf{m}_4=[\bar{1}11]$, respectively. For both these anisotropic cases, parameters are $\omega_i=6$, $\rho_i=0.9$, $\gamma_0=1$, and $\beta=0.001$.
}
    \label{fig:figure6}
    \end{figure}

As a last example we demonstrate the possibility to simulate coarsening of bulk nanoporous materials, as e.g. considered with $\gamma(\hat{\mathbf{n}}) = 1$ in \citen{Geslin2019,BeckAndrewsetal_arXiv_2020}, in a strongly anisotropic setting. We highlight a design project for the Digital Archive of Mathematical Models (DAMM)~\cite{damm}. Within computational domains, defined by the Wulff shapes, a classical Cahn-Hilliard model with random perturbations of $u = \frac{1}{2}$ as initial condition, is solved. After spinodal decomposition the solution is used as initial condition for the anisotropic variational DDCH model, which is solved for several time steps similarly to Fig.~\ref{fig:figure5}. The isosurface $u = \frac{1}{2}$ is extracted, postprocessed and 3D printed. Fig. \ref{fig:figure6} shows three examples of the resulting objects, a sphere (isotropic) and a cube and a tetrahedra with rounded corners and edges.   

\section{Conclusions}
\label{sec:conclusions}

In this paper we have extended the variational DDCH model for 
isotropic surface diffusion\cite{SalvalaglioDCH} to the anisotropic case. We consider weak and strong anisotropies.
The first case only requires the energy with a singular restriction function to be multiplied by a surface energy $\gamma(\hat{\mathbf{n}})$, as in \citen{Torabi2009}. The second case requires an additional regularisation
by a Willmore functional, again following \citen{Torabi2009}. A direct connection of the resulting models with the
well known, and well used, non-variational \textit{RRV model}, \citen{Ratz2006}, as possible in the isotropic case, cannot be achieved 
for anisotropic surface energies. However, numerical comparisons show the similarities of both approaches in terms of
approximation properties to the sharp interface limit. We omit the asymptotic analysis of the anisotropic variation DDCH model. Formal convergence results to anisotropic
surface diffusion can be obtained by combining the analysis in \citen{SalvalaglioDCH,Ratz2006,Torabi2009}.

The slightly more complex variational DDCH models, if compared with the non-variational \textit{RRV model}, \citen{Ratz2006}, not only gives
energy dissipation and a mathematical foundation for numerical analysis, it also provides better positivity preserving properties.
Besides the superior approximations properties, this positivity preserving property enables the wide applicability of the model in \citen{Ratz2006} for large-scale applications. Our numerical results indicate that this property is further improved in the variational DDCH models. This is beneficial for several reasons. For a numerical point of view it allows to optimize criteria based on thresholds or changes of $u$ as, e.g., for refining spatial discretization or adaptive timestepping. Also, these models are often coupled with additional equations in the bulk phases, which are characterized by the limiting values of $u$, e.g., elasticity\cite{Ratz2006,Salvalaglio2018} or compositions\cite{Bergamaschini2020}. In these cases, any improvement on positivity preservation, leads to increased accuracy. 

However, besides these advantages, which will further foster to application of these models in materials science, as for the isotropic setting, several questions related to existence, uniqueness,
and regularity of solutions, as well as to the positivity preserving property, remain open and need to be addressed. \\

    \section*{Acknowledgements}
This research was partially funded by the EU H2020 FET-OPEN project microSPIRE (ID: 766955) and by the EU H2020 FET-OPEN project NARCISO (ID: 828890). We gratefully acknowledge the computing time granted by the John von Neumann Institute for Computing (NIC) and provided on the supercomputer JURECA at J\"ulich Supercomputing Centre (JSC), within the Project no. HDR06, and by the Information Services and High Performance Computing (ZIH) at the Technische Universit\"at Dresden (TUD). SMW acknowledges generous financial support from the US National Science Foundation through grants NSF-DMS 1719854 and NSF-DMS 2012634.

    \section*{Conflicts of Interest Statement}
The authors certify that they have no affiliations with or involvement in any organization or entity with any financial or non-financial interest in the subject matter discussed in this manuscript.

%\bibliographystyle{alpha}
%\bibliography{references}

\end{document}